\documentclass{article}
\usepackage{amsfonts}
\usepackage{amssymb}
\usepackage{latexsym}
\def\be{\begin{equation}}
\def\ee{\end{equation}}
\def\bea{\begin{eqnarray}}
\def\eea{\end{eqnarray}}
\def\qed{\hfill\mbox{$\Box$}\medskip}
\def\Qed{\hfill\rule{1ex}{1ex}\medskip}
\def\PP{{\cal P}}
\def\notdiv{\nmid}

\newcommand{\Z}{\mathbf{Z}}
\newcommand{\0}{^{\phantom{0}}}
\begin{document}
\begin{Large}
\centerline{\bf Every even number greater than $454$}
\centerline{\bf is the sum of seven cubes}
\end{Large}
\vspace*{2ex}
\centerline{Noam D.~Elkies}
\vspace*{2ex}
\centerline{September, 2010}

\vspace*{5ex}

\begin{quote}
{\bf Abstract.}  It is conjectured that every integer $N>454$
is the sum of seven nonnegative cubes.  We prove the conjecture
when $N \equiv 2 \bmod 4$.  This result, together with a recent proof
for $4|N$, shows that the conjecture is true for all even~$N$.
\end{quote}

{\large\bf 1 Introduction}

Linnik [1943] 
showed that every sufficiently large natural number~$N$\/
is the sum of at most seven positive cubes.  It has long been known that
the set of $N>0$ without such a representation contains
\bea
  & \{15,\; 22,\; 23,\; 50,\; 114,\; 167,\; 175,\; 186,\; 212,\;
   \qquad\qquad
\label{eq:exc7} \\
  & \qquad\qquad
   231,\; 238,\; 239,\; 303,\; 364,\; 420,\; 428,\; 454\},
\nonumber
\eea
and conjectured that (\ref{eq:exc7}) is the full set of exceptions.
See the first section of~\cite{Ram07} for this history of this part of
Waring's problem.   The introduction of~\cite{BE} gave a briefer account;
an even more abbreviated summary follows.

The first effective upper bound on the largest exception
was $\exp(\exp(13.94))$, obtained in \cite{McCurley};
the upper bound now stands at $\exp(524)$
by the analytic sieve argument of~\cite{Ram07}.
While any effective upper bound reduces the problem to a finite computation,
$\exp(524)$ is still much too large to reach by such a computation
in practice.  But it is still useful to have a good lower bound
on any exceptional $N$\/ outside the known set (\ref{eq:exc7}),
because known techniques for constructing seven-cube representations
typically require $N$\/ to be somewhat large.  We shall use the bound
$N \geq 2.5 \cdot 10^{26}$, proved in~\cite{BRZ}.  The largest such
bound reported is $\exp(78.7) > 1.5 \cdot 10^{34}$ \cite[pages 433--434]{DHL},
which is still very far below $\exp(524) > 3 \cdot 10^{224}$.

A different tack is to prove the conjecture under some
congruence condition on~$N$, in the hope that eventually every $N$\/
might be covered by one such result.  We know of two such theorems.
The first \cite{BRZ} shows that $N$\/ is the sum of seven cubes if
$N \equiv 0$ or $\pm1 \bmod 9$ and $N$\/ is an invertible cubic residue
mod~$37$, and also that $37$ could be replaced by
some larger primes congruent to $1 \bmod 3$.
The second \cite{BE} proves that $N$\/ is the sum of seven cubes if $4|N$\/
as long as $N$\/ is outside the set (\ref{eq:exc7}) of known exceptions.

Here we prove this result for $2 \| N$, which together with \cite{BE}
establishes it whenever $2|N$:

{\bf Theorem.}
{\em If $N$\/ is an even positive integer not in
\be
\{22,\;50,\;114,\;186,\;212,\;238,\;364,\;420,\;428,\;454\}
\label{eq:exc7_even}
\ee
then $N$\/ is the sum of seven nonnegative cubes.}

The proof adapts a technique used in earlier work on
Waring's problem for cubes, dating back to the initial paper 
\cite{Wieferich} on sums of nine cubes, and including \cite{BRZ,BE}.
We need only consider $N \equiv 2 \bmod 4$
since the case $4|N$\/ was proved in~\cite{BE}.
We use two new ingredients:

$\bullet$ We take coefficients $(a_1,a_2,a_3) = (1,2,5)$
in \hbox{$Q = \sum_{i=1}^3 a_i\0 X_i^2$}.
As with other quadratic forms used in such constructions,
this $Q$\/ is diagonal and unique in its genus,
but one doesn't expect a factor of~$5$ in $a_1 a_2 a_3$
when there is no condition on $N \bmod 5$.
(Factors of~$3$ occur in~\cite{BE}, but that is not surprising
because of the special behavior of cubes
in the \hbox{$p$-adic} integers $\Z_p$ for $p=3$:
they are all congruent to $0$ or $\pm 1 \bmod 9$.)

$\bullet$ Due to the structure of the set of positive integers
not represented by~$Q$, we must restrict the auxiliary parameter~$p$
to a residue class modulo~$300$.
This modulus is beyond the range of the tables of \cite{RamRum}.
Extending these tables to primes in congruence classes mod~$300$
would require a large computation with Dirichlet \hbox{$L$-functions}.
Instead we replace the prime~$p$ by a product~$P$\/ of distinct primes
each congruent to $5 \bmod 6$.
This retains the key property that every residue class has a cube root,
while giving enough flexibility to reduce the use of~\cite{RamRum}
to primes in two arithmetic progressions mod~$12$.
This refinement also lets us dispense with the factor $\beta$ of~\cite{BE},
since it can be included among the prime factors of~$P$.
It also streamlines or completes several other seven-cube constructions
of this kind in~\cite{BE2}, where we prove that $N$\/ is the sum of
seven cubes if $N \equiv 0 \bmod 9$, $\pm 1 \bmod 18$, or $\pm 2 \bmod 9$.
Note that the first two of these,
together with the results of \cite{BE} and the present paper,
properly contain the theorem of \cite{BRZ} by
removing the additional hypothesis modulo $37$ or a larger prime.

It is noted in \cite{BE} that the construction actually produces
a representation $N = \sum_{i=0}^6 x_i^3$ with each $x_i$ positive,
not just nonnegative, once $N$\/ is large enough, with $N=2408$ probably
being the last exception.  The same is true here; indeed one can easily
adjust the proof to produce a representation with $\max_i x_i / \min_i x_i$
uniformly bounded: it is enough to replace the bounds $1618$ and $1786$
on $N/P^3$ by $1618+\delta$ and $1786-\delta$, and and to change $x_0\0$ to
$x'_0 = x_0\0 + 6P$ if $x_0\0 < \delta^2 P$,
for sufficiently small $\delta > 0$.
The same can be done for the results in \cite{BRZ} and~\cite{BE}.

The rest of this paper is organized as follows.  In the next section
give the formulas (\ref{eq:458},\ref{eq:N458},\ref{eq:Q125})
that represent $N$\/ as a sum of seven cubes given suitable $P,Q$.
The following section derives conditions on~$P$\/ that guarantee
that the criteria on~$Q$\/ can be satisfied (the new analysis here
is in Lemma~1).  Finally we prove (Lemma~2) that such $P$\/
can be found if $N>10^{20}$; this together with the bound
$2.5 \cdot 10^{26}$ of~\cite{BRZ} completes the proof of the theorem.

\vspace*{2ex}

{\large\bf 2 From $N$\/ to $Q$}

Given $N \equiv 2 \bmod 4$ with $N$\/ large enough, namely $N > 10^{20}$,
we shall construct a representation of~$N$\/ as the sum of seven cubes.
We may assume $5^3 \notdiv N$, because we may write $N = 5^{3e} N_0$
with $5^3 \notdiv N_0$, and use a seven-cube decomposition
$N_0 = \sum_{i=0}^6 x_i^3$ to write $N$\/ as the sum of the cubes of
$5^e x_i$.  If $N_0$ is in the exceptional set (\ref{eq:exc7}),
but $e>0$ (so $N \neq N_0$), we use a decomposition of $5^3 N_0$
to the same effect; indeed for each of these $17$ values of $N_0$
a direct computation shows that five cubes suffice, as do
seven positive cubes (this is used to represent $N$\/ as a sum of
seven positive cubes, and is contained in the computation reported
in~\cite{BE} that suggests all $N>2408$ have such a representation).

We start from the usual six-cube identity (see e.g.\ \cite[Lemma~1]{BE}),
taking $(c_1,c_2,c_3) = (4,5,8)$ to find

\vspace*{2ex}
\centerline{
$
  (4P + X_1)^3 + (4P - X_1)^3
 + (5P + X_2)^3 + (5P - X_2)^3
 + (8P + X_3)^3 + (8P - X_3)^3
$
}
\vspace*{-2.5ex}
\be
{} = 1402 P^3 + 6P Q_1,
\label{eq:458}
\ee
where $Q_1 := 4 X_1^2 + 5 X_2^2 + 8 X_3^2$.  We deduce that
if $P$\/ is a positive integer, and $X_1, X_2, X_3$ are integers such that
$|X_1| < 4P$, $|X_2| < 5P$, and $|X_3| < 8P$,
then for any positive integer~$x_0$ we have a representation of
\be
N = x_0^3 +  1402 P^3 + 6P Q_1
\label{eq:N458}
\ee
as a sum of seven positive cubes.

In our setting, $N \equiv 2 \bmod 4$, so we must take $x_0$ even.
We shall require that $P$\/ be odd.  Then $1402 P^3 \equiv 2 \bmod 4$,
so $6P Q_1$ is a multiple of~$4$, whence $Q_1$ is even.
Therefore $2 | X_2$, from which $4 | Q_1$.
We thus have $Q_1 = 4 Q$\/ where $Q$\/ is the diagonal quadratic form
\be
Q := X_1^2 + 2 X_3^2 + 5 (X_2/2)^2.
\label{eq:Q125}
\ee
Dickson proved that this~$Q$\/ represents all nonnegative integers
except those of the form $5^{2k} (25 n \pm 10)$,
see \cite[p.69, Theorem~IX]{Dickson}.\footnote{
  The theorem statement follows the proof on pages 67--69;
  the notation $G$\/ for this form is on page~63.
  }
We shall show that for $N > 10^{20}$,
we can choose $P$\/ such that there is an $x_0$ that makes $Q$\/
a positive integer, not congruent to~$0$ or $\pm 10 \bmod 25$,
with $Q < 16 P^2$.  Then $4 Q < 4^3 P^2$, so for each $i=1,2,3$ we have
\be
4 X_i^2 \leq c_i X_i^2 \leq Q_1 = 4 Q < 4^3 P^2,
\label{eq:Xineq}
\ee
so $|X_i| < 4 P \leq c_i P$, and the $X_i$ satisfy
the inequalities that make each term positive in the seven-cube
representation of~$N$\/ obtained from formulas (\ref{eq:458},\ref{eq:N458}).

\pagebreak

{\large\bf 3 The conditions on~$P$}

Taking $Q_1 = 4 Q$\/ in (\ref{eq:N458}) and solving for~$Q$\/ yields
\be
Q = \frac{N - x_0^3 - 1402 P^3}{24P}.
\label{eq:Q}
\ee
We noted already that $x_0$ must be even and $P$\/ odd.  We now see
that $P$\/ must also be chosen so that
\be
P \equiv \frac{N}{2} \bmod 4
\label{eq:Pmod4}
\ee
to make the numerator of~(\ref{eq:Q}) divisible by~$8$.
To make $Q$\/ integral, it remains to choose $x_0$ so that
$3P | N - x_0^3 - 1402 P^3$.  To that end, we require that
$P$\/ be a product of distinct primes each congruent to $5 \bmod 6$.
Then every residue mod~$6P$\/ has a cube root.  We choose for $x_0$
the smallest positive solution of $x_0^3 \equiv N - 1402P^3 \bmod 6P$.
Then $x_0$ is an even number in $(0,6P]$, and $Q$\/ is an integer in
\be
\left(
  \frac{N - (1402+6^3) P^3}{24P}, \; \frac{N - 1402 P^3}{24P}
\right].
\label{eq:Q_bound}
\ee
Therefore $Q > 0$ provided $N/P^3 > 1402 + 6^3 = 1618$,
and $Q < 16 P^2$ provided $N/P^3 < 1402 + 24 \cdot 16 = 1786$.
Thus we seek $P$\/ in the interval $(A N^{1/3}, B N^{1/3})$ where
$A = 1786^{-1/3}$ and $B = 1618^{-1/3}$,
with $B / A = (893 / 809)^{1/3} > 1.033$.
(These estimates correspond to Proposition~1 in~\cite{BE}.)

To the condition (\ref{eq:Pmod4}) on $P\bmod 4$, we next add
a condition mod~$25$ to assure that $Q$\/ is congruent to neither
$0$ nor $\pm 10 \bmod 25$, from which it will follow that
$Q$\/ is represented by the diagonal quadratic form with coefficients $1,2,5$.

{\bf Lemma 1.} {\em
For any $N \not\equiv 0 \bmod 5^3$ there exist at least two
nonzero choices of $b \in \Z/25\Z$ such that if $P \equiv b \bmod 25$
then (\ref{eq:Q}) yields a value of $Q$\/ not congruent to
$0$ or $\pm 10 \bmod 25$.
}

(Note that since $b$\/ is not the zero residue
the congruence $P \equiv b \bmod 25$ cannot force $P$\/ to have
a repeated prime factor.)

{\em Proof}\/:
If $5|N$\/ we shall choose either $b = \pm 5$ or $b = \pm 10$.
Since $5|P$, also $5|x_0$, and then the numerator of (\ref{eq:Q}) is
$$
N - x_0^3 - 1402 P^3 \equiv N \bmod 5^3.
$$
In particular, if $25 \notdiv N$\/
then either $b = \pm 5$ or $b = \pm 10$ works, because
the numerator is not a multiple of~$25$,
so its quotient by $24P$\/ is not a multiple of~$5$,
and thus lies outside the forbidden congruence classes mod~$25$.
If $25 | N$\/ then we choose $b = \pm 5$ if $N \equiv 25$ or $100 \bmod 5^3$,
and $b = \pm 10$ if $N \equiv 50$ or $75 \bmod 5^3$
(recall that $5^3 \notdiv N$\/).
Then $Q \equiv N / (24P) \equiv \pm 5 \bmod 25$, so again
$Q$\/ is outside the exceptional set for $X_1^2 + 2 X_2^2 + 5 X_3^2$.

Finally, if $5 \notdiv N$, we first choose $b \bmod 5$ so that
$1402 b^3 \equiv N \bmod 5$, then choose $b \bmod 25$ so that
\be
\frac{N - 1402 b^3}{24 b} \equiv \pm 5 \bmod 25.
\label{eq:b_cond}
\ee
The equation for $b \bmod 5$ has a solution because every integer is
a cube mod~$5$; since $N$\/ is not a multiple of~$5$, neither is~$b$.
The choice of $b \bmod 5$ guarantees that $5|Q$\/ if and only if
$5|x_0$, in which case $Q \equiv (N - 1402 b^3) / (24 b) \bmod 5^3$.
Thus (\ref{eq:b_cond}) means that if $5|Q$\/ then $Q \equiv \pm 5 \bmod 25$.
We claim that each sign arises for a unique choice of $b \bmod 25$,
which we call $b_\pm$.  We tabulate $b_+$ and $b_-$ for each of the
$20$ possible residues of $N \bmod 25$:

\vspace*{1ex}

{\small
\centerline{
\setlength{\arraycolsep}{.3em}
$
\begin{array}{c|cccccccccccccccccccc}
N \bmod 25 \, & \,
 1 &  2 &  3 &  4 &  6 &  7 &  8 &  9 & 11 & 12 & 13 & 14 & 16 & 17 & 18 & 19 & 21 & 22 & 23 & 24
\\ \hline
b_+ &
 2 & 21 &  9 & 18 & 22 &  1 & 14 & 13 &
17 &  6 & 19 &  8 & 12 & 11 & 24 &  3 &
 7 & 16 &  4 & 23
\\ \hline
b_- &
 7 &  6 & 24 & 13 &  2 & 11 &  4 &  8 &
22 & 16 &  9 &  3 & 17 & 21 & 14 & 23 &
12 &  1 & 19 & 18
\end{array}
$
}
}
\vspace*{1ex}
\centerline{Table 1: $b_+$ and $b_-$ for each $N \bmod 25$}

This completes the proof of Lemma~1.\qed

{\em Remark}\/: The existence and uniqueness of $b_\pm$ in each case
can be understood as follows.  Fix an arbitrary $b_0$ such that
$N - 1402 b_0^3 \equiv 0 \bmod 5$, and let $Q_0 = (N - 1402 b_0^3) / (24 b_0)$.
Then $b = b_0 + 5 \beta$ yields $Q \equiv Q_0 + 5 Q'_0 \beta \bmod 25$ where
$$
Q'_0 = \frac{-N - 2 \cdot 1402 b_0^3}{24 b_0^2}
$$
is the $5$-adic derivative of~$Q$\/ at~$b_0$ (so the linear approximation
$Q_0 + 5 Q'_0 \beta$ is within $O(5^2)$ of the correct value).
But a function of the form $(A-B b^3)/b$ cannot vanish at some~$b$\/
together with its derivative $(-A-2B^3)/b^2$ modulo~$5$
(or indeed mod~$p$ for any prime $p \neq 3$).
Since $Q_0 \equiv 0 \bmod 5$, then, $Q'_0 \not \equiv 0 \bmod 5$,
so $Q_0 + 5 Q'_0 \beta \bmod 25$ runs over all lifts of~$Q_0$ to $\Z/25\Z$
as $\beta$ varies mod~$5$.

\vspace*{2ex}

{\large\bf 4 The existence of~$P$}

It remains to prove:

{\bf Lemma 2.}  {\em
Suppose $N \equiv 2 \bmod 4$ with $125 \notdiv N$, and $N>10^{20}$.
Choose $b$ according to Lemma~1.
Then there exists $P \in (AN^{1/3}, BN^{1/3})$
that is a product of distinct primes congruent to $5 \bmod 6$
with $P \equiv b \bmod 25$ and $P \equiv (N/2) \bmod 4$.
}

{\em Proof}\/:  Suppose first that $N > 10^{26}$.
We then choose $P$\/ as follows.  Let $\PP$\/ be the set of
squarefree integers each of whose prime factors is congruent to $5 \bmod 6$.
Choose a finite subset $\PP_0 \subset \PP$\/ such that:
\begin{itemize}
\item All elements of $\PP_0$ are congruent mod~$4$, say to $r_0$;
\item $\PP_0$ contains a representative of each nonzero class mod~$25$;
\item $\max(\PP_0)/\min(\PP_0) = 1 + \epsilon_0\0$ with $\epsilon_0\0$ small;
\item $\max(\PP_0)$ is as small as possible given $\epsilon_0\0$.
\end{itemize}
We shall take $P = P_0 p$ with $P_0 \in \PP_0$
and $p$ a prime greater than $\max(\PP_0)$ that is congruent to
$5 \bmod 6$ and to $r_0 (N/2) \bmod 4$.  We choose $p$ first,
and then select $P_0$ so that $P \equiv b \bmod 25$
with $b$\/ depending on $N \bmod 25$ according to Lemma~1.
Then $P$\/ is a product of distinct primes each congruent to $5 \bmod 6$,
and is in a suitable class mod~$4$ and mod~$25$ to guarantee the success
of our construction as long as $1618 < N/P^3 < 1786$.  This last condition,
in turn, will be satisfied provided that
\be
p \in (p_{\min},p_{\max}) := \left(
  \frac{(N/1786)^{1/3}}{\PP_{\min}},
  \frac{(N/1618)^{1/3}}{\PP_{\max}}
\right),
\label{eq:p_gap}
\ee
an interval whose endpoints' ratio is
\be
\frac{p_{\max}}{p_{\min}} = \frac{ (1786/1618)^{1/3} }{\PP_{\max}/\PP_{\min}}
>
\frac{1.033}{1+\epsilon_0\0}.
\label{eq:p_gap_ratio}
\ee

A quick computer search finds the choice
\be
\PP_0 = \PP \cap (1 + 4\Z) \cap [26141,26669],
\label{eq:P0}
\ee
with $\epsilon_0\0 = 528/26141 < 0.0202$ and
$1.033 / (1+\epsilon_0\0) > 1.0125$.
We tabulate the $38$ elements $P_0 \in \PP_0$,
sorted by the remainder $\overline{\! P}_0$ of $P_0 \bmod 25$,
together with the prime factorization of those $P_0 \in \PP_0$
that are not prime:

\vspace*{1ex}

\centerline{\small
$
\begin{array}{r|l}
\overline{\! P}_0 \! & P_0 \\ \hline
1 & 26401 = 17 \cdot 1553; \ 26501
\\
2 & 26177; \ 26477 = 11 \cdot 29 \cdot 83
\\
3 & 26153; \ 26653 = 11 \cdot 2423
\\
4 & 26329 = 113 \cdot 233
\\
5 & 26305 = 5 \cdot 5261
\\
6 & 26281 = 41 \cdot 641
\\
7 & 26357
\\
8 & 26633
\\
9 & 26309; \ 26609 = 11 \cdot 41 \cdot 59
\\
10 & 26185 = 5 \cdot 5237; \ 26485 = 5 \cdot 5297
\\
11 & 26261; \ 26461 = 47 \cdot 563; \ 26561
\\
12 & 26237
\end{array}
\qquad
\begin{array}{r|l}
\overline{\! P}_0 \! & P_0 \\ \hline
13 & 26513
\\
14 & 26189; \ 26389 = 11 \cdot 2399; \ 26489
\\
15 & 26365 = 5 \cdot 5273; \ 26665 = 5 \cdot 5333
\\
16 & 26141 \ [\min]
\\
17 & 26417
\\
18 & 26393
\\
19 & 26669 \ [\max]
\\
20 & 26345 = 5 \cdot 11 \cdot 479; \ 26545 = 5 \cdot 5309
\\
21 & 26321; \ 26521 = 11 \cdot 2411
\\
22 & 26297; \ 26597
\\
23 & 26473 = 23 \cdot 1151; \ 26573
\\
24 & 26249
\end{array}
$
}
\vspace*{1ex}
\centerline{Table 2:
   $P_0 \in [26141, 26669]$ for each nonzero $\overline{\! P}_0 \bmod 25$
}

\vspace*{1ex}

Suppose first that $p_{\min} > 10^{10}$.  We then take $k=12$ in
\cite[Theorem~1]{RamRum}, finding that for each $l \in (\Z/12\Z)^*$
there exists a prime $p \equiv l \bmod 12$ such that
$1 < p/p_{\min} < (1+\epsilon_{12}\0) / (1-\epsilon_{12}\0)$,
with $\epsilon_{12}\0 < 0.003$ according to \cite[Table~1, p.~419]{RamRum}.
This proves Lemma~2 for all
\be
N > 1786 (10^{10} \PP_{\max})^3 = 3.38767\ldots \cdot 10^{46}.
\label{eq:N46}
\ee
If $\PP_{\max} < p_{\min} < 10^{10}$, we apply the algorithm used in
the Sublemma of \cite[Lemma 2]{BE} for primes congruent to $l \bmod 12$
for $l=5$ and $l=11$.  We quickly find that for any
$p_{\min} \in (26669,10^{11})$ the interval $(p_{\min}, 1.006 p_{\min})$
contains at least one prime $p \equiv 5 \bmod 12$
and at least one prime $p \equiv 11 \bmod 12$.
(Indeed in these two arithmetic progressions mod~$12$
the largest ratio between consecutive primes past $\PP_{\max}$ is
$35381/35201 < 1.00512$ for $l=5$ and $45491/45263 < 1.00504$ for $l=11$.)
This extends the range of our proof from (\ref{eq:N46}) down to
\be
N > 1786 \PP_{\max}^6 = 6.42572\ldots \cdot 10^{29}.
\label{eq:N29}
\ee
[Note that this is already sufficient to prove our theorem,
using the bound $10^{34}$ of \cite{Ram07, DHL}.]

Finally, to bring the bound on~$N$\/ from (\ref{eq:N29}) down to
the claimed $10^{20}$, we apply the same algorithm directly to $P$\/
in each of the $48$ odd residue classes mod~$100$ that is not 
a multiple of~$25$, searching not for primes congruent to $5 \bmod 6$
but for elements of~$\PP$.  We soon find that for any
$P_{\min} \in ( (10^{20}/1786)^{1/3}, \PP_{\max} )$
the interval $(P_{\min}, (1+\epsilon_0\0) P_{\min})$
contains at least one $P \in \PP$\/ in each of those
$48$ classes.  This completes the proof of Lemma~2 and of our theorem.\Qed

{\em Remarks}\/:  Indeed in each of those $48$ classes the ratio between
two consecutive $P \in \PP$\/ never gets as large as $1.015$.
We could have dispensed with the step from (\ref{eq:N46}) to (\ref{eq:N29})
by extending the last computation past $\PP_{\max}$ to $10^{10}$;
the required calculation would be much larger (more residue classes, and
prime factorization rather than just primality testing), but still feasible.
On the other hand we could not use the same argument to reduce
the bound on~$N$\/ to the $10^{18}$ used in~\cite{BE},
because the ratio $1.0389+$ between the consecutive elements
$92437 = 23 \cdot 4019$ and $96037 = 137 \cdot 701$ of
$\PP \cap (100\Z + 37)$ is too large.  We could probably use the choice between
$b_+$ and $b_-$ (see Table~1) to extend the range of our construction
down to $10^{18}$, but not much lower.

\vspace*{2ex}

{\large\bf Acknowledgements}

I thank Kent Boklan and Ali Assarpour for a computer file listing
of all $102$ regular diagonal forms $Q = \sum_{i=1}^3 a_i\0 X_i^2$
together with the arithmetic progressions not represented by
each of these forms~$Q$.\footnote{
  The list was obtained by Jones in his doctoral thesis \cite{Jones},
  which though unpublished can be found on Jagy's quadratic-forms
  webpages, see the Bibliography.  This list was first published in
  \cite{JonesPall}, together with proofs but without the lists of
  excluded progressions.  That information appears in
  \cite[\S58, Table 5, pages 112--113]{Dickson39}; see also page~111
  for the abbreviations $A$, $B$, \ldots,~$N$\/ for the sets appearing
  in this table.  Jagy also reproduces this table online, see
  {\sf http://zakuski.math.utsa.edu/\~{}kap/Forms/Dickson\_Diagonal\_1939.pdf}~.
  }

I also thank Boklan, as well as Nathan Kaplan, for a careful reading of
an earlier draft of this paper.

This work is partly supported by the
National Science Foundation under grant DMS-501029.

\begin{small}

\textsc{Department of Mathematics, Harvard University,
Cambridge, MA 02138, U.S.A.} (\textsf{elkies@math.harvard.edu})

\end{small}

\end{document}